\documentclass[11pt,a4paper]{amsart}
\usepackage{amssymb,amsmath}
\numberwithin{equation}{section}
\vspace{8cm}
\newtheorem{theorem}{Theorem}
\newtheorem{lemma}[theorem]{Lemma}
\newtheorem{prop}[theorem]{Proposition}
\newtheorem{cor}[theorem]{Corollary}

\newtheorem{oss}[theorem]{Remark}

\numberwithin{theorem}{section}
\newtheorem*{theorem*}{Theorem}
\newcommand{\Pk}{\mathcal{P}^+_k}

\newcommand{\R}{{\mathbb R}}

\newcommand{\refe}[1]{{(\ref{#1})}}

\newcommand{\noi}{\noindent}

\newcommand{\e}{\'e\ }

\begin{document}
\parindent=0pt

\title[Entire subsolutions...]{Entire subsolutions of fully nonlinear\\ degenerate elliptic equations}

\author[I. Capuzzo Dolcetta, F. Leoni, A. Vitolo]
{Italo Capuzzo Dolcetta, Fabiana Leoni, Antonio Vitolo}

\address{ ICD and FL: Dipartimento di Matematica\newline
\indent Sapienza Universit\`a  di Roma}
\email{capuzzo@mat.uniroma1.it}
\email{leoni@mat.uniroma1.it}
\address{AV: Dipartimento di  Matematica\newline
Universit\`a di Salerno}
\email{vitolo@unisa.it}

\begin{abstract}
We prove existence and non existence results for fully nonlinear degenerate elliptic inequalities, by showing that the classical Keller--Osserman condition on the zero order term is a necessary and sufficient condition for the existence of entire sub solutions.

 \end{abstract}\maketitle

\section{Introduction}\label{INTRO}

Consider the semilinear equation
\begin{equation}\label{Brezis}
\Delta u = |u|^{\gamma-1}u+g(x)
\end{equation}
with $\gamma>1$ and $g(x)\geq \varepsilon>0$ bounded and continuous.
We know from Brezis \cite{B} that this equation is uniquely solvable in $\mathbb R^n$ and, by known regularity and comparison results, there is a solution $u<0$ in $\R^n$. Therefore, the function $v=-u$ is a solution to
$$
\Delta v=|v|^\gamma -g(x)\, .
$$
\noindent Consider now the equation
\begin{equation}\label{Osserman-eq}
\Delta u = |u|^\gamma+g(x)\, ,
\end{equation}
and observe that if $u$ is a solution of the above then $u$ solves also
\begin{equation}\label{Osserman}
\Delta u \ge f(u)
\end{equation}
where
$$
f(t):=\begin{cases}
t^{\gamma}+\varepsilon & \hbox{ if }\ \ t\ge 0\\
\varepsilon &  \hbox{ if } \ \ t< 0
\end{cases}
$$
is a positive, non decreasing and  continuously differentiable function such that
$$
\int_0^{+\infty}\left(\int_0^{t}f(s)ds\right)^{-\frac12}dt = \int_0^{+\infty}\left((\gamma+1)^{-1}t^{\gamma+1} +\varepsilon t\right)^{-\frac12}dt < +\infty \,\, .
$$

Therefore, the Keller--Osserman condition
\begin{equation}\label{KO}
\int_0^{+\infty}\left(\int_0^{t}f(s)ds\right)^{-\frac12}dt = +\infty
\end{equation}
is not verified and from well-known results by Keller \cite{K} and Osserman \cite{O}, we deduce that  inequality (\ref{Osserman}), and therefore also equation \refe{Osserman-eq}, cannot have entire solutions. \\

We are interested here in investigating the validity of this type of results for fully nonlinear  degenerate elliptic inequalities. More precisely, we consider viscosity solutions of the  partial differential inequality
\begin{equation}\label{EQ1}
F(D^2u) \ge  f(u)\quad\hbox{in}\;\;\R^n
\end{equation}
where $F$ is a second order  degenerate elliptic operator in the sense of Crandall, Ishii, Lions \cite{CIL} and $f(u)$ is a positive, non decreasing zero order term. \\
Of particular interest are those mapping $F:\mathcal{S}_n\to \R$, where $\mathcal{S}_n$ is the space of $n\times n$ real symmetric matrices, which are functions of the eigenvalues.
In our model cases, $F$ will be either  the  elliptic operator $\Pk$  defined for any  $X\in \mathcal{S}_n$ and a positive integer $1\leq k\leq n$ as
\begin{equation}\label{OP1}
\Pk (X) = \mu_{n-k+1}(X) +\ldots +\mu_n(X) =\sup_{W\in G(k,n)} \mathrm{Trace}_W (X) 
\end{equation}
$\mu_1(X) \leq \mu_2 (X) \leq \ldots \leq \mu_n(X)$ being the ordered eigenvalues of the matrix $X$ and $G(k,n)$ being the Grassmanian of $k$--dimensional subspaces of $\R^n$,  or  the degenerate maximal Pucci operator defined by
\begin{equation}\label{OP2}
\mathcal{M}^+_{0,1} (X) = \sum_{\mu_i>0} \mu_i(X)\, = \sup_{A\in \mathcal{S}_n: A\leq I_n} \mathrm{Trace}\, (AX)
\end{equation}

\noi Let us point out however that these two operators do not belong to the class of degenerate Hessian operators investigated by Neil S. Trudinger and collaborators, see in this respect the remarkable results about new maximum principles and regularity in e.g. \cite{TRUD1,TRUD2,NewMaxPrin}.

The Pucci extremal operators have been extensively  studied  by  Caffarelli, Cabr\e in the uniformly elliptic case, see  \cite{CAFCA}. Let us recall here that the operator \refe{OP2} is maximal not only in the class of linear operators, but it also bounds from above all degenerate elliptic operators vanishing at $O$. In particular, for any $1\leq k\leq n$  and for all $X\in \mathcal{S}_n$ one has 
$$
\Pk (X)\leq \mathcal{M}^+_{0,1} (X)\,.
$$
As for the operators
$\Pk$,  we refer to the recent works of  Caffarelli, Li, Nirenberg \cite{CLN, CLN2}, Harvey, Lawson \cite{HL1,HL2,HL3}, see also Amendola, Galise, Vitolo \cite{AGV}, and  the references therein. We just point out here that  such degenerate operators arise in several frameworks, e.g. the  geometric problem of mean curvature evolution of manifolds with co-dimension greater than one, as in Ambrosio, Soner \cite{AS}, as well as   the PDE approach to the convex envelope problem,  see Oberman, Silvestre \cite{OS}.

After the above mentioned classical results  in \cite{K}, \cite{O}, \cite{B}  about entire solutions of the semi linear equation (\ref{Osserman}), many extensions have been obtained for different operators and more general zero order terms. In particular, for divergence form principal parts let us  recall  the results of Boccardo, Gallouet, Vazquez \cite{BGV1, BGV2}, Leoni \cite{L} and Leoni, Pellacci \cite{LP}, D'Ambrosio, Mitidieri \cite{DM}. In the fully nonlinear framework, analogous results have been more recently obtained by Esteban, Felmer, Quaas \cite{EFQ}, Diaz \cite{D}  and Galise, Vitolo \cite{GV}, and by  Bao, Ji \cite{BJ}, Bao, Ji, Li \cite{BJL}, Jin, Li, Xu \cite{JLX} for Hessian equations, involving the $k$-th elementary symmetric function of the eigenvalues $\mu_1(D^2u), \ldots ,\mu_n(D^2u)$.\\
In these papers, existence, uniqueness and comparison results are given for the equation
 $$
 F(D^2u) = f(u)-g(x)\, 
 $$
 under local integrability assumptions on the datum $g$ and by assuming the zero order term $f$ to be of absorbing type. For example,  $f:\R \to \R$ is odd, continuous, increasing, convex for $t\geq 0$ and satisfying  the growth condition for $t\to\infty$
 $$
\int^{+\infty}\left(\int_0^{t}f(s)ds\right)^{-\frac12}dt <+\infty
$$
In the present paper, we complement the already established results by considering the different case in which $f$ is bounded from below, say positive, and non decreasing, and $F$ is degenerate elliptic.

Our main results are the following ones:
 \begin{theorem}\label{main}
Let $1\leq k\leq n$ and $f:\mathbb R \to \mathbb R$ be  positive, continuous and  non decreasing. Then the inequality 
\begin{equation}\label{eqpk}
\Pk (D^2u)\geq f(u)
\end{equation}
has an entire viscosity solution $u\in C(\R^n)$  if and only if $f$ satisfies the Keller-Osserman condition \refe{KO}.
\end{theorem}

 \begin{theorem}\label{mainbis}
Let $f:\mathbb R \to \mathbb R$ be  positive, continuous and  strictly increasing. Then the inequality 
\begin{equation}\label{eqpk2}
\mathcal{M}^+_{0,1} (D^2u)\geq f(u)
\end{equation}
has an entire viscosity solution $u\in C(\R^n)$  if and only if $f$ satisfies the Keller-Osserman condition \refe{KO}.
\end{theorem}

The proof of both theorems is based, as in the semi linear case, on  a comparison argument with radial symmetric functions obtained as solutions of an associated ODE.\\
Remarkably, the comparison principle  works also in the present cases where degeneracy eventually occurs both in the principal part and in the zero order term. \\ 
Let us observe that, by the maximality of the operator $\mathcal{M}^+_{0,1}$, Theorem \ref{mainbis} gives a necessary condition for the existence of entire viscosity solutions of 
$$
F(x,D^2u)\geq f(u)\, ,
$$
for any continuous operator $F:\R^n\times \mathcal{S}_n\to \R$ satisfying $F(x,O)=0$ and the ellipticity condition
\begin{equation}\label{ell}
0\leq F(x,X+Y)-F(x,X)\leq \mathrm{Trace}(Y)
\end{equation}
for all $x\in \R^n$ and $X, \ Y\in \mathcal{S}_n$ with $Y\geq O$.
\\ Moreover, Theorem \ref{mainbis} combined with the above mentioned results of \cite{D} provides a necessary and a sufficient condition for the existence and uniqueness of an entire viscosity solution of the non homogeneous equation
$$
\mathcal{M}^+_{0,1} (D^2u)=f(u)-g(x)
$$
with a bounded and continuous datum $g$, see Corollary \ref{g}.

 \noi Let us finally mention that  the arguments used to prove  Theorem \ref{main} and Theorem \ref{mainbis} can be adapted to more general partial differential inequalities involving first order terms, and  we refer in this respect to our work in progress \cite{CDLV}.

\section{On the ODE \;$\varphi''(r)+\frac{c-1}r\,\varphi'(r)=f(\varphi(r))$}

\noindent  In order to obtain existence/non existence results for viscosity solutions of inequalities (\ref{eqpk}), (\ref{eqpk2}), let us first  investigate the associated 
 second order ODE 
\begin{equation}\label{ODE}
\varphi''(r)+\frac{c-1}r\,\varphi'(r)=f(\varphi(r))\,, \ \ r \ge 0
\end{equation}
for a given positive constant $c>0$ and equipped with the initial  condition
\begin{equation}\label{I'}
\varphi'(0)=0\,.
\end{equation}
By a solution $\varphi \in C^2\left( [0,R)\right) $ of \refe{ODE}, \refe{I'}  we mean a function $\varphi \in C^2\left( (0,R)\right) $ with $0 < R \le +\infty$, which is continuous in $[0,R)$,  twice differentiable in $(0,R)$  and such that
$$ 0=\varphi'(0)=\lim_{r \to 0^+}\varphi'(r)\,,\, \varphi''(0)=\lim_{r \to 0^+}\varphi''(r)=\lim_{r\to 0^+}\frac{\varphi'(r)}{r}\neq \infty\,.$$

\begin{lemma}\label{convex-increasing} Let  $c>0$ and $f:\R\to \R$ be  continuous, non negative  and non decreasing. Then, every solution $\varphi $ of problem \refe{ODE}, \refe{I'} in some interval $[0,R)$ is non decreasing and convex. 
\end{lemma}
\noindent{\bf Proof}.  Writing \refe{ODE} in the form
\begin{equation}\label{equivalent-eqn}
(r^{c-1}\varphi')'= f(\varphi(r))r^{c-1}\, 
\end{equation}
we see that  $r^{c-1}\varphi'(r)$ is  non decreasing since $f$ is non-negative. Hence, $\varphi'(r)\geq 0$ for $r\geq0$. \\
Next, integrating (\ref{equivalent-eqn}) between $0$ and $s$, using the assumption $\varphi'(0)=0$ and the monotonicity of $f\circ \varphi$ we have
$$
s^{c-1}\varphi'(s)= \int_0^s(r^{c-1}\varphi')'dr= \int_0^sf(\varphi(r))r^{c-1}dr \le \frac{s^c}c\,f(\varphi(s))\,,
$$
from which 
\begin{equation}\label{eigenvalue-comparison}
\frac{\varphi'(s)}{s} \le \frac{f(\varphi(s))}{c}\,. 
\end{equation}
Using this information in equation \refe{ODE} we get
\begin{equation}\label{last}
\varphi''(s) = f(\varphi(s))- \frac{c-1}s\,\varphi'(s)\ge \frac{\varphi'(s)}{s}\geq 0\, ,
\end{equation} 
showing that   $\varphi$ is convex. \hfill$\Box$

\begin{oss}\label{rem-ineq} {\rm If, in addition, $f$ is strictly positive, then every solution $\varphi$  of \refe{ODE}, \refe{I'} will be accordingly strictly increasing  and strictly convex. Moreover, observe that if $c \geq1$ then for $s \in [0,R)$ 
\begin{equation}\label{ineqs}
\frac{f(\varphi(s))}{c} \le \varphi''(s) \le f(\varphi(s))
\end{equation}
Indeed, the left-hand inequality follows from \refe{eigenvalue-comparison} inserted into \refe{last} while the right-hand inequality is obtained from equation \refe{ODE} observing that $c \geq1$ and  we have just proved that $\varphi'\ge 0$.}
\end{oss}

\noindent The existence of  solutions of equation \refe{ODE} follows from classical ODE theory with continuous data. As for the maximal interval of existence we have the following result due to Osserman \cite{O}, whose proof is included here for the reader's convenience.

\begin{lemma}\label{maximality} Let  $c\geq1$. If $f:\R \to \R$ is continuous, non negative  and  non decreasing, then every maximal solution of \refe{ODE}, \refe{I'} is globally defined in $[0,+\infty)$ if and only if $f$ satisfies \refe{KO}.
\end{lemma}
\noi {\bf  Proof}.  If $f\equiv 0$, then \refe{KO} is trivially fulfilled. On the other hand, in this case solutions of \refe{ODE}, \refe{I'} are necessarily constants, since $c\geq1$.  Assume now that $f$ does not vanish identically, and let $\varphi :[0,R)\to \R$ be a non constant maximal solution of problem \refe{ODE}, \refe{I'}. Then, there exists $r_0\in [0,R)$ such that $\varphi'(r)>0$ for $r\geq r_0$ and, by \refe{eigenvalue-comparison}, $f(\varphi(r))>0$ for $r\geq r_0$. Multiplying (\ref{ineqs}) by $\varphi'(s)$,   and then integrating between $r_0$ and $r$ we get
$$
\frac2{ c}\int_{r_0}^rf(\varphi(s))\varphi'(s)ds + \left( \varphi'(r_0)\right)^2\le \left( \varphi'(r)\right) ^2\le 2 \int_{r_0}^rf(\varphi(s))\varphi'(s)ds + \left( \varphi'(r_0)\right)^2\, .
$$
Since  $\varphi$ is a $C^1$-diffeomorphism between $(r_0,R)$ and $ (\varphi(r_0), \varphi(R))$,   it follows that
$$
\frac2{ c}\int_{\varphi(r_0)}^{\varphi(r)}f(t)dt +\left( \varphi'(r_0)\right)^2\le \left(\varphi'(r)\right)^2\le 2 \int_{\varphi(r_0)}^{\varphi(r)}f(t)dt +\left( \varphi'(r_0)\right)^2\, ,
$$
that is
$$
\left(2 \int_{\varphi(r_0)}^{\varphi}f(t)dt+ \left( \varphi'(r_0)\right)^2\right)^{-\frac12}\le r'(\varphi)\le \left(\frac2{c}\int_{\varphi(r_0)}^{\varphi}f(t)dt +\left( \varphi'(r_0)\right)^2)\right)^{-\frac12}\, .
$$
Integrating between $\varphi(r_0)$ and $\varphi(R)$ yields
$$
\int_{\varphi(r_0)}^{\varphi(R)}\frac{d\varphi}{\sqrt{2 \int_{\varphi(r_0)}^{\varphi}f(t)dt+ \left( \varphi'(r_0)\right)^2}} \le R-r_0\le \int_{\varphi(r_0)}^{\varphi(R)}\frac{d\varphi}{\sqrt{2/c  \int_{\varphi(r_0)}^{\varphi}f(t)dt+ \left( \varphi'(r_0)\right)^2}}
$$

Therefore, if $R=+\infty$, then the right hand side integral is infinite positive and necessarily $\varphi (R)=+\infty$ and the Keller-Osserman condition (\ref{KO}) is satisfied.\\
Assume conversely  that (\ref{KO}) holds, and suppose by contradiction  that $R<+\infty$. Since $[0,R)$ is the maximal interval of existence of the monotonically non decreasing  solution $\varphi(r)$, we  have $\varphi(r) \to +\infty$ as $r \to R^-$, so that the first of above inequalities yields a contradiction with (\ref{KO}).
\hfill$\Box$  

\begin{oss}\label{f>01} {\rm We observe that if $f$ is as in Lemma \ref{maximality} and satisfies
\begin{equation}\label{nKO}
\int^{+\infty}\left( \int_0^tf(s)\, ds\right)^{-1/2}dt<+\infty\, ,
\end{equation}
then problem \refe{ODE}, \refe{I'} can have in general both  maximal solutions globally existing in $[0,+\infty)$, namely  constant solutions,  and maximal solutions defined only on a bounded interval $[0,R)$. But,  if we assume $f$ to be strictly positive in $\R$, then conditions \refe{ODE}, \refe{I'} do not allow for constant solutions, and the above proof shows that  either \emph{all} maximal solutions  are global or \emph{all}  maximal solutions are defined on a bounded subset of $[0,+\infty)$ according to condition \refe{KO} is satisfied or not. In particular, if $f>0$ and satisfies \refe{nKO}, then every maximal solution $\varphi$ cannot be defined beyond $[0, R)$ with $R$ satisfying}
$$
R\leq \int_{\varphi(0)}^{+\infty} \sqrt{\frac{c}{2\int_{\varphi(0)}^\varphi f(t)dt}}d\varphi\; .
$$
\end{oss}

\section{Fully nonlinear degenerate elliptic inequalities}

\noindent By using Lemma \ref{convex-increasing}, we are in position to show now that a classical solution of equation 
\begin{equation}\label{Pucci-soln}
\Pk (D^2\Phi) = f(\Phi)\quad \hbox{ in } B_R
\end{equation}
 can be obtained from a solution $\varphi \in C^2\left([0,R)\right)$ of problem (\ref{ODE}), (\ref{I'}) with $c=k$
by setting $$\Phi(x)=\varphi(|x|), \ \ |x|<R\,.$$ 

\begin{lemma}\label{L1} Let $1\leq k\leq n$,
 $f:\mathbb R \to \mathbb R$ be  non negative, non decreasing and continuous  and $\varphi \in C^2\left([0,R)\right)$ be a solution of problem \refe{ODE}, \refe{I'} with $c=k$.\\ Then $\Phi(x)=\varphi(|x|)$ is a  classical  solution of equation \refe{Pucci-soln}. 
\end{lemma}
\noindent {\bf Proof}.  We notice that 
$$
D^2\Phi (x)=\left\{
\begin{array}{ll}
\varphi''(0)\, I_n & \hbox{ if } x=0\\[2ex]
\frac{\varphi'(|x|)}{|x|}\, I_n + \left( \varphi''(|x|)-\frac{\varphi'(|x|)}{|x|}\right) \frac{x}{|x|}\otimes \frac{x}{|x|} & \hbox{ if } x\neq 0
\end{array} \right. 
$$
Hence, it is easy to check that  $\Phi \in C^2(B_R)$, and that the eigenvalues of $D^2\Phi (x)$ are $\varphi''(0)$, with multiplicity $n$ if $x=0$, and $\varphi'' (|x|)$, which is simple, and $\frac{\varphi'(|x|)}{|x|}$ with multiplicity $n-1$ for $x\neq 0$. 
 We then have
$$
\Pk (D^2\Phi (0))=k\, \varphi''(0)=f(\varphi (0))=f(\Phi(0))
$$
and, by inequality \refe{last} in  Lemma \ref{convex-increasing},
$$
\Pk (D^2\Phi(x))=\varphi''(|x|)+\frac{k-1}{|x|}\,\varphi'(|x|)=f(\varphi(|x|))=f(\Phi(x))\qquad \hbox{for } x\neq 0
$$
as well, so that $\Phi$ is a classical solution of equation (\ref{Pucci-soln}). 
\hfill$\Box$

\noindent In the next result we establish  a form of comparison principle between merely continuous viscosity subsolutions and smooth supersolutions which holds true even in the currently considered degenerate case, see also at this purpose \cite{G}, \cite{AGV}.

\begin{prop}\label{comp} Assume $f:\R \to \R$ continuous and nondecreasing, let $u\in C(B_R)$ and $\Phi\in C^2(B_R)$ be, respectively, a viscosity subsolution and  a classical supersolution of \refe{Pucci-soln}.\\ If
$$\limsup_{|x| \to R^-}\,\left(u(x)-\Phi(x)\right)\le 0$$ 
then $u(x) \le \Phi(x)$ for all $x \in B_R$\,.  
\end{prop}
\noindent {\bf Proof}. By contradiction, suppose there is some point $x \in B_R$ where $u(x)>\Phi(x)$. Hence, taking $\varepsilon>0$ small enough, we have  that the set $\Omega:\equiv\{x \in B_R, \ u(x)-\Phi(x)>\varepsilon\}$ is non-empty  and that $\overline \Omega \subset B_R$.\\
Set $v(x) = u(x)-\Phi(x)$ in $B_R$. Since $u\in C(B_R)$ is a viscosity subsolution  and $\Phi$ is a classical supersolution, one has
$$
\Pk (D^2v) \ge \Pk (D^2u)-\Pk (D^2\Phi)\ge f(u)-f(\Phi)
$$
in the viscosity sense in $B_R$.
 Since $f$ is non decreasing, it follows that $\Pk (D^2v)\ge 0$ in $\Omega$. Moreover, $v>\varepsilon$ in $\Omega$ $v=\varepsilon$ on $\partial \Omega$. Hence, there exists a concave paraboloid $\Psi (x)$ touching $v$ from above at some point  $x_0\in \Omega$, a contradiction to the inequality $\Pk (D^2 \Psi (x_0))\geq 0$.
 \hfill$\Box$

\begin{oss}\label{OSS1} {\rm  As a comparison function $\Phi$ in Proposition \ref{comp}, one can take  $\Phi(x)=\varphi(|x|)$ where $\varphi\in C^2\left([0,R)\right)$ is  any convex non decreasing solution  of }
$$
\begin{cases}
\varphi''+\frac{k-1}{r}\,\varphi'\le f(\varphi) \ \  \hbox{in }\ \ [0,R)\\
\varphi'(0)=0
\end{cases}
$$
 \end{oss}

\begin{oss}\label{Puccio} {\rm Let us observe that if we strengthen the assumption on $f$ by requiring its  strict monotonicity, then  the above proof works as well for the degenerate maximal Pucci operator in (\ref{OP2}) yielding the validity of the comparison principle in Proposition \ref{comp}  for this strongly degenerate elliptic operator.}
 \end{oss}

Combining Lemma \ref{L1} and Proposition \ref{comp} with the maximality result of Lemma \ref{maximality}, we can now show that  Keller-Osserman condition (\ref{KO}) is a necessary and sufficient condition for the existence of entire solutions of the differential inequalities (\ref{eqpk}) and (\ref{eqpk2}).

\noi {\bf Proof of Theorem \ref{main}}.  Assume that $(\ref{eqpk})$ has a viscosity solution $u\in C(\R^n)$ and let  $\varphi\in C^2\left( [0,R)\right)$ be a maximal solution of \refe{ODE}, \refe{I'} satisfying  the extra initial condition $\varphi(0) <u(0)$.  We claim that $R=+\infty$. If, on the contrary, $R < +\infty$ then $\varphi(r) \to +\infty$ as $r \to R^-$ and  $\Phi(x)=\varphi(|x|)$ blows up on the boundary $\partial B_R$. Hence, $u(x) \le \Phi(x)$ in $B_R$ by Proposition \ref{comp}, a contradiction to $u(0)> \varphi(0)$.\\  Therefore, the maximal interval of existence of $\varphi$ is $[0,+\infty)$ and, by Lemma \ref{maximality} and Remark \ref{f>01},  condition (\ref{KO}) is satisfied.\\
Conversely, suppose that the Keller-Osserman condition \refe{KO} holds true and let $\varphi$ be a maximal solution of \refe{ODE}, \refe{I'}. Again, Lemma \ref{maximality} implies that $\varphi$ is globally defined on $[0,+\infty)$ and, by Lemma \ref{L1}, that $u(x)=\varphi(|x|)$ is an entire classical solution of \refe{eqpk}. 
\hfill$\Box$

\noi {\bf Proof of Theorem \ref{mainbis}}. Use Remark \ref{Puccio} and proceed exactly as in the proof above. \hfill$\Box$

\indent {\rm Let us discuss now the more general case where the strict positivity condition on $f$ in Theorems \ref{main}, \ref{mainbis} is relaxed to $f\geq 0$. In this case, there exists $t_0\in \R$ such that $f(t)\equiv 0$ for $t\leq t_0$ and $f(t)>0$ for $t>t_0$. Then, inequality \refe{eqpk} has, of course, entire constant solutions $u(x)\equiv c$ for any $c\leq t_0$ and  one may ask about existence of non-constant entire solutions. \\
Looking at the  proof of Lemma \ref{maximality}, we see that if $f$ satisfies the Keller--Osserman condition \refe{KO}, then the ODE problem \refe{ODE}, \refe{I'} does have indeed non-constant global solutions $\varphi$, namely those solutions satisfying the initial condition $\varphi(0)>t_0$. By Lemma \ref{L1}, any such non-constant global solution $\varphi$ generates an entire non-constant  solution $u$ of \refe{eqpk}.\\
 On the other hand, the same argument used in  the proof of Theorem \ref{main} shows that if there exists an entire solution $u$ of \refe{eqpk}  such that $u(x_0)>t_0$ at some point $x_0\in \R^n$, then $f$ must satisfy \refe{KO}. \\
Therefore, in order to show that, in the present case, \refe{KO} is a necessary and sufficient condition for the existence of non-constant entire solutions of \refe{eqpk}, one has to prove the validity of a Liouville  type theorem for the operator $\Pk$, stating  the non existence of non-constant bounded from above solutions of
\begin{equation}\label{Lio}
\Pk (D^2u)\geq 0\quad \hbox {in }\ \R^n\, .
\end{equation}
For $k=n\leq2$, the classical Liouville theorem for subharmonic functions applies and we get the conclusion. On the contrary, if $n\geq 3$ or $n=2$ and $k=1$, then inequality \refe{Lio} admits non-constant  solutions bounded from above, namely any smooth radial  function $u(x)=\varphi (|x|)$ with $\varphi$ bounded and increasing and with $u$ subharmonic if $n\geq 3$.\\
Therefore, in these cases, \refe{KO} is a sufficient but not a necessary condition for the existence of non-constant entire solutions of \refe{eqpk}.
\\
\indent Let us finally recall, see \cite{CL}, that a  Liouville theorem holds true for the uniformly elliptic Pucci's inf--operator
$$
\mathcal{M}^-_{\lambda, \Lambda} (X)=\lambda\, \sum_{\mu_i>0}\mu_i(X) +\Lambda\, \sum_{\mu_i<0}\mu_i(X)\, 
$$
with ellipticity constants $\Lambda \geq \lambda >0$, provided the space dimension $n$ satisfies the restriction $n\leq 1+\frac{\Lambda}{\lambda}$.\\ Observing that all the other arguments used in the proof of Theorem \ref{main} can be applied also for the operator $\mathcal{M}^-_{\lambda, \Lambda}$, from the previous discussion we deduce the validity of the following statement:

\begin{prop} Let $\Lambda\geq \lambda>0$ and  $f:\R \to \R$  continuous, non decreasing and non negative.  If $n\leq 1+\frac{\Lambda}{\lambda}$, then there exist non constant solutions of
$$
\mathcal{M}^-_{\lambda, \Lambda}(D^2u)\geq f(u)\quad \hbox{in }\ \R^n
$$
if and only if $f$ satisfies \refe{KO}.
\end{prop}
\noindent By means of the next result and recalling Remark \ref{f>01} we can enlarge the class of functions $f$ for which the Keller--Osserman condition \refe{KO} is a necessary condition for the existence of entire solutions of \refe{eqpk}:

\begin{cor}\label{t0}
Let $1\leq k\leq n$ and $f:\R \to \R$ be positive, continuous, non decreasing for $t\geq t_0$ and satisfying \refe{nKO}.
 Then, there does not exist any entire viscosity solution of inequality \refe{eqpk}.
 \end{cor}
\noi {\bf Proof}.   The function
$$
\tilde{f}(t)=\left\{
\begin{array}{ll}
f(t) & \hbox{ if } t\geq t_0\\[1ex]
\min_{[t,t_0]} f(s) & \hbox{ if } t<t_0
\end{array} \right.
$$
satisfies all the assumptions in Theorem \ref{main} as well as \refe{nKO}, so that there does not exist any entire viscosity solution $u$ of
$$
\Pk (D^2u)\geq \tilde{f}(u)\, .
$$
Since $f(t)\geq \tilde{f}(t)$ for any $t\in \R$, it follows that no entire viscosity solution of \refe{eqpk} can exists as well.
\hfill$\Box$

\noindent Combining the results we have obtained with  previously known results for equations having strictly increasing absorbing zero order terms, we finally  deduce the following existence/non existence statement for viscosity solutions of the non homogeneous equation 
\begin{equation}\label{eqg}
\mathcal{M}^+_{0,1} (D^2u) =f(u) - g(x)\, ,\qquad x\in \R^n\, .
\end{equation}

\noi The typical case covered by the next result is $f(u)=\exp u\,.$
\begin{cor}\label{g}
Let $f:\R \to \R$  be continuous, strictly increasing, convex, bounded from below and $g:\R^n\to \R$ a bounded, continuous function. Assume also that
\begin{equation}\label{nKO2}
\int^{+\infty} \frac{dt}{\sqrt{\int_0^t \left( f(s)-\inf_\R f\right) ds}}< +\infty\, .
\end{equation}
 Then\\[1ex]
(i) if $\sup_{\R^n} g\leq \inf_\R f$, then \refe{eqg} does not have any viscosity solution,  \\[1ex]
(ii) if $\inf_{\R^n} g > \inf_\R f$, then \refe{eqg} has a unique bounded viscosity solution .
\end{cor}

\noi {\bf Proof}. Statement (i) immediately follows from Theorem \ref{mainbis} with $f$ replaced by $$\tilde{f} (t)=f(t)-\inf_\R f\,.$$ As far as (ii) is concerned,  we observe that, by  assumption, there exists $t_0\in \R$ such that $g(x)\geq f(t_0)$ for all $x\in \R^n$. Let us consider the function
$$
\tilde{f}(t)=\left\{
\begin{array}{ll}
f(t+t_0)-f(t_0) & \hbox { if } t\geq 0\\[1ex]
-\tilde{f}(-t) & \hbox { if } t< 0\, ,
\end{array}
\right.
$$
which is continuous, odd, increasing and convex for $t\geq 0$. A convexity argument shows that 
$\tilde{f}$ satisfies, for all $t\in \R$ and $h\geq 0$, the inequality
\begin{equation}\label{beta}
\tilde{f}(t+h)-\tilde{f}(t)\geq 2\, \tilde{f}\left( \frac{h}{2}\right)\, .
\end{equation}
Using \refe{beta} and \refe{nKO2}, 
we can apply results in Diaz \cite{D} to deduce the existence of a unique bounded viscosity solution $v\in C(\R^n)$ of
\begin{equation}\label{mmtilde}
\mathcal{M}^+_{0,1} (D^2v)=\tilde{f}(v)+ f(t_0)-g(x)\, ,\qquad x\in \R^n\, .
\end{equation}
Moreover, by comparison and the assumptions made on $g$, we have $v\geq 0$, so that $u(x)=v(x)+t_0$ is a bounded viscosity solution of \refe{eqg}. \\ We finally observe that, if $u$ and $v$ are two bounded viscosity solution of \refe{eqg}, then, for $t_0=\min \{ \inf_{\R^n}u,\ \inf_{\R^n}v\}$, both $u-t_0$ and $v-t_0$ solve \refe{mmtilde}. By the uniqueness proved in \cite{D}, we then conclude that $u\equiv v$.
\hfill$\Box$

\begin{oss}
{\rm It is easy to check that Corollary \ref{g} holds true for any principal part of the form $F(x,D^2u)$, with $F:\R^n\times \mathcal{S}_n\to \R$ continuous, satisfying $F(x,O)=0$ and the ellipticity condition \refe{ell}.
}
\end{oss}


\bigskip
\bigskip
\begin{thebibliography}{30}
\bibitem{AS} L. Ambrosio, H. M. Soner,  Level set approach to mean curvature flow in arbitrary codimension, {\it J. Differential Geom. } \textbf{43} (4),  693--737  (1996)
\bibitem{AGV} M.E. Amendola, G. Galise and A. Vitolo, Riesz capacity, maximum principle and removable sets of fully nonlinear second order elliptic operators, preprint
\bibitem{BJ} J. Bao, X. Ji, Necessary and sufficient conditions on solvability for Hessian inequalities, {\it Proc. Amer. Math. Soc.} \textbf{138}, 175--188 (2010)
\bibitem{BJL} J. Bao, X. Ji, H. Li, Existence and nonexistence theorem for entire subsolutions of $k$--Yamabe type equations, {\it J. Differential Equations} \textbf{253}, 2140--2160 (2012)
\bibitem{BGV1} L. Boccardo, T. Gallouet, J.L. Vazquez, Nonlinear elliptic equations in $\R^N$ without growth restriction on the data, {\it J. Differential Eq.} \textbf{105} (2), 334--363 (1993)
\bibitem{BGV2} L. Boccardo, T. Gallouet, J.L. Vazquez, Solutions of nonlinear parabolic equations without growth restrictions on the data, {\it Electr. J.Differential Eq.} \textbf{2001} (60), 1--20 (2001)
\bibitem{B} H. Brezis, Semilinear equations in $\mathbb R^n$ without conditions at infinity, {\it Appl. Math. Optim.} \textbf {12}, 271--282 (1984)
\bibitem{CAFCA}  L.A. Caffarelli, X. Cabr\e, {\it Fully Nonlinear Elliptic Equations},  American Mathematical Society Colloquium Publications \textbf{43}  (1995)
\bibitem{CLN} L.A. Caffarelli, Y.Y .Li and L. Nirenberg, Some remarks on singular solutions of nonlinear elliptic equations. I, {\it J. Fixed Point Theory Appl.} \textbf{5}, 353--395 (2009)
\bibitem{CLN2} L.A. Caffarelli, Y.Y. Li and L. Nirenberg, Some remarks on singular solutions of nonlinear elliptic
equations. III: viscosity solutions, including parabolic operators, {\it Comm. Pure Appl. Math.} doi:
10.1002/cpa.21412 (2012)
\bibitem{CDLV} I. Capuzzo Dolcetta, F. Leoni, A. Vitolo, in preparation
\bibitem {CIL} M. G. Crandall, H. Ishii, P. L. Lions, User's guide to
viscosity solutions of second order
partial differential equations, {\it Bulletin of the American Mathematical
Society} \textbf {27} (1), 1--67
(1992)
\bibitem{TRUD2} N. Chauduri, N.S. Trudinger, An Alexsandrov type theorem for k convex functions, {\it Bull. Austral. Math. Soc.} \textbf{71}, 305-314 (2005)
\bibitem{CL} A. Cutr\`\i , F. Leoni, On the Liouville property for
fully nonlinear equations,  {\it Ann.Inst. H. Poincar\e,
Analyse Nonlineaire} \textbf{17} (2), 219--245 (2000)
\bibitem{DM} L. D'Ambrosio, E. Mitidieri, 
 A priori estimates, positivity results, and nonexistence theorems for quasilinear degenerate elliptic inequalities, {\it Advances in Mathematics} \textbf{224},  967-1020 (2010)
\bibitem{D} G. Diaz, A note on the Liouville method applied to elliptic eventually degenerate fully nonlinear equations governed by the Pucci operators and the Keller--Osserman condition, {\it Math. Ann.} \textbf{353}, 145--159 (2012)
\bibitem{EFQ} M.J. Esteban, P.L. Felmer, A. Quaas, Super-linear elliptic equations for fully nonlinear operators without growth restrictions for the data, {\it Proc. R. Soc. Edinb.} \textbf{53}, 125--141 (2010)
\bibitem{G} G. Galise, Maximum principles, entire solutions and removable singularities of fully nonlinear second order equations, {\it Ph.D. Thesis}, Universit\`a\  di Salerno a.a. 2011/2012
\bibitem{GV} G. Galise, A. Vitolo, Viscosity solutions of uniformly elliptic equations without boundary and growth conditions at infinity, {\it Int. J. Differ. Equ.} \textbf{2011}  (2011)
\bibitem{TRUD1} P. Guan, N.S. Trudinger, X-J. Wang, On the Dirichlet problem for degenerate Monge-Amp\`ere equations, {\it Acta Math.} \textbf{182}, 87-104 (1999)
\bibitem{HL1} R. Harvey,  B. Lawson Jr., Dirichlet duality and the nonlinear Dirichlet problem, {\it Comm. Pure
Appl. Math.} \textbf{62}, 396 443 (2009)
\bibitem{HL2} R.Harvey and B.Lawson Jr., Plurisubharmonicity in a general geometric context, {\it Geometry and Analysis} \textbf{1}, 363--401 (2010)
\bibitem{HL3} R.Harvey and B.Lawson Jr., Dirichlet duality and the nonlinear Dirichlet problem on Riemannian manifolds,
{\it J. Diff. Geom.} \textbf{88} , 395--482 (2011)
\bibitem{HL4} R.Harvey and B.Lawson Jr., Existence, uniqueness and removable singularities for nonlinear partial differential equations in geometry, to appear in  {\it Surveys in Geometry}, ArXiv:1303.1117
\bibitem{JLX} Q. Jin, Y.Y. Li, H. Xu, Nonexistence of positive solutions for some fully nonlinear elliptic equations, {\it Methods Appl. Anal.} \textbf{12}, 441--449 (2005)
\bibitem{K} J.B. Keller, On solutions of $\Delta u=f(u)$, {\it Comm. Pure Appl. Math.} \textbf{10}, 503--510 (1957)
\bibitem{NewMaxPrin} H.Kuo and N.Trudinger, New maximum principles for linear elliptic equations, {\it Indiana University Mathematics Journal} \textbf{56}, no. 5, 2439--2452 (2007) 
\bibitem{L} F. Leoni, Nonlinear elliptic equations in $\R^N$ with \lq\lq absorbing" zero order terms, {\it Adv. Differ. Equ.} \textbf{5}, 681--722 (2000)
\bibitem{LP} F. Leoni, B. Pellacci, Local estimates and global existence for strongly nonlinear parabolic equations with locally integrable data, {\it J. Evol. Equ.} \textbf{6},  113�-144 (2006)
\bibitem{OS} A. Oberman, L. Silvestre, The Dirichlet problem for the convex envelope, {\it Trans. Amer. Math. Soc. } \textbf{363} (11), 5871--5886  (2011)
\bibitem{O} R. Osserman,  On the inequality $\Delta u \ge f(u)$, {\it Pacific J. Math.} \textbf{7} ,  1141--1147 (1957)

\end{thebibliography}
\end{document}